\documentclass{article}
\usepackage{amsmath}
\usepackage{amsfonts}
\usepackage{amssymb}
\usepackage{graphicx}
\newtheorem{theorem}{Theorem} [section]

\newtheorem{proposition}[theorem]{Proposition}

\newenvironment{proof}[1][Proof]{\noindent\textbf{#1.} }{\ \rule{0.5em}{0.5em}}
\begin{document}

\author{Vadim E. Levit\\Department of Computer Science and Mathematics\\Ariel University Center of Samaria, Israel\\levitv@ariel.ac.il
\and Eugen Mandrescu\\Department of Computer Science\\Holon Institute of Technology, Israel\\eugen\_m@hit.ac.il}
\title{Greedoids on Vertex Sets of Unicycle Graphs}
\date{}
\maketitle

\begin{abstract}
A \textit{maximum stable set }in a graph $G$ is a stable set of maximum size.
$S$ is a \textit{local maximum stable set} of $G$, and we write $S\in\Psi(G)$,
if $S$ is a maximum stable set of the subgraph spanned by $S\cup N(S)$, where
$N(S)$ is the neighborhood of $S$. $G$ is a \textit{unicycle graph }if it owns
only one cycle. In \cite{LevMan2} we have shown that the family $\Psi(T)$ of a
forest $T$ forms a greedoid on its vertex set. Bipartite, triangle-free, and
well-covered graphs $G$ whose $\Psi(G)$ form greedoids were analyzed in
\cite{LevMan4,LevMan6,LevMan08b}, respectively.

In this paper we characterize the unicycle graphs whose families of
local maximum stable sets form greedoids.


\end{abstract}

\section{Introduction}

\thispagestyle{empty}

Throughout this paper $G=(V,E)$ is a simple (i.e., a finite, undirected,
loopless and without multiple edges) graph with vertex set $V=V(G)$ and edge
set $E=E(G)$. If $X\subset V$, then $G[X]$ is the subgraph of $G$ induced by
$X$, and by $G-W$ we mean the subgraph $G[V-W]$, where $W\subset V(G)$. The
graph $G$ is \textit{unicycle }if it owns only one cycle. The
\textit{neighborhood} of a vertex $v\in V$ is the set $N(v)=\{w:w\in V,vw\in
E\}$. If $\left\vert N(v)\right\vert =1$, then $v$ is a \textit{pendant
vertex}. We denote the \textit{neighborhood} of the set $A\subset V$ by
$N_{G}(A)=\{v\in V-A:N(v)\cap A\neq\emptyset\}$ and its \textit{closed
neighborhood} by $N_{G}[A]=A\cup N_{G}(A)$, or shortly, $N(A)$ and $N[A]$, if
there is no ambiguity. 

By $K_{n},C_{n}$ we mean the complete graph
on $n\geq1$ vertices, and the chordless cycle on $n\geq3$ vertices, respectively.

A \textit{stable} set in $G$ is a set of pairwise non-adjacent vertices. A
stable set of maximum size will be referred to as a \textit{maximum stable
set} of $G$, and the \textit{stability number }of $G$, denoted by $\alpha(G)
$, is the cardinality of a maximum stable set in $G$. By $\Omega(G)$ we denote
the family of all maximum stable sets of the graph $G$.

\begin{figure}[h]
\setlength{\unitlength}{1.0cm} \begin{picture}(5,1.3)\thicklines
\multiput(4.5,0)(1,0){5}{\circle*{0.29}}
\multiput(5.5,1)(1,0){3}{\circle*{0.29}}
\put(4.5,0){\line(1,0){4}}
\put(5.5,0){\line(0,1){1}}
\put(6.5,1){\line(1,0){1}}
\put(6.5,0){\line(0,1){1}}
\put(7.5,1){\line(1,-1){1}}
\put(4.5,0.35){\makebox(0,0){$a$}}
\put(5.3,0.35){\makebox(0,0){$b$}}
\put(6.2,0.35){\makebox(0,0){$c$}}
\put(7.5,0.35){\makebox(0,0){$d$}}
\put(6.2,1.2){\makebox(0,0){$g$}}
\put(7.8,1.2){\makebox(0,0){$h$}}
\put(8.5,0.35){\makebox(0,0){$e$}}
\put(5.3,1.2){\makebox(0,0){$f$}}
\end{picture}
\caption{A graph having {various local maximum stable sets}.}%
\label{fig10}%
\end{figure}
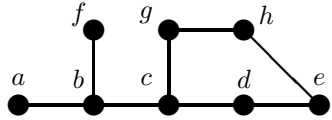

A set $A\subseteq V(G)$ is a \textit{local maximum stable set} of $G$ if $A$
is a maximum stable set in the subgraph induced by $N[A]$, i.e., $A\in
\Omega(G[N[A]])$, \cite{LevMan2}. Let $\Psi(G)$ stand for the family of all local
maximum stable sets of $G$. For instance, any set $S$ consisting of only
pendant vertices belongs to $\Psi(G)$, while the converse is not generally
true; e.g., the set $\{e,g\}\in\Psi(G)$ contains no pendant vertex, where $G$
is the graph in Figure \ref{fig10}.

Clearly, not any stable set of a graph $G$ is included in some maximum stable
set of $G$. For example, there is no $S\in\Omega(G)$ such that
$\{b,d,h\}\subset S$, where $G$ is the graph presented in Figure \ref{fig10}.
In \cite{NemhTro}, Nemhauser and Trotter Jr. showed that every local maximum
stable set of a graph can be enlarged to one of its maximum stable sets.

A \textit{matching} in a graph $G=(V,E)$ is a set of edges $M\subseteq E$ such
that no two edges of $M$ share a common vertex. We denote the size of a
\textit{maximum matching} (a matching of maximum cardinality) by $\mu(G)$. Let
us recall that $G$ is a \textit{K\"{o}nig-Egerv\'{a}ry graph }provided
$\alpha(G)+\mu(G)=\left\vert V(G)\right\vert $, \cite{Dem,Ster}. 
It is known that every bipartite graph is a K\"{o}nig-Egerv\'{a}ry graph \cite{eger,koen}.

A \textit{greedoid}, \cite{BjZiegler,KorLovSch}, is a pair
$(V,\mathcal{F})$, where $\mathcal{F}\subseteq2^{V}$ is a non-empty set system
satisfying the following conditions:

\setlength {\parindent}{0.0cm}\textit{accessibility}: for every non-empty
$X\in\mathcal{F}$, there is an $x\in X$ such that $X-\{x\}\in\mathcal{F}%
$;\setlength
{\parindent}{3.45ex}

\setlength {\parindent}{0.0cm}\textit{exchange}: for $X,Y\in\mathcal{F}%
,\left\vert X\right\vert =\left\vert Y\right\vert +1$, there is an $x\in X-Y$
such that $Y\cup\{x\}\in\mathcal{F}$.\setlength
{\parindent}{3.45ex}\newline

The following theorem shows that it is enough to prove that $\Psi(G)$
satisfies the accessibility property, in order to validate that $\Psi(G)$ is a greedoid.

\begin{theorem}
\label{th1}\cite{LevMan08c} If the family $\Psi(G)$ satisfies the
accessibility property, then it satisfies the exchange property as well.
\end{theorem}

Clearly, $\Omega(G)\subseteq\Psi(G)$ holds for any graph $G$.

If $S\in\Psi(G)$, $\left\vert S\right\vert =k\geq2$, then sometimes there
exists a chain
\[
\{x_{1}\}\subset\{x_{1},x_{2}\}\subset...\subset\{x_{1},...,x_{k-1}%
\}\subset\{x_{1},...,x_{k-1},x_{k}\}=S
\]
such that $\{x_{1},x_{2},...,x_{j}\}\in\Psi(G)$, for all $j\in\{1,...,k-1\}$;
such a chain is called an \textit{accessibility chain }for $S $,
\cite{LevMan4}. It is evident that $x_{1}$ must be a simplicial vertex, i.e., a vertex
whose closed neighborhood induces a complete graph in $G$ (in particular, any
pendant vertex is also simplicial). 

\begin{figure}[h]
\setlength{\unitlength}{1.0cm} 
\begin{picture}(5,1.1)\thicklines
\multiput(2,0)(1,0){4}{\circle*{0.29}}
\multiput(3,1)(1,0){3}{\circle*{0.29}}
\put(2,0){\line(1,0){3}}
\put(5,0){\line(0,1){1}}
\put(3,1){\line(1,0){1}}
\put(3,0){\line(0,1){1}}
\put(4,0){\line(0,1){1}}
\put(2,0.3){\makebox(0,0){$a$}}
\put(2.7,1){\makebox(0,0){$b$}}
\put(3.7,0.3){\makebox(0,0){$c$}}
\put(5.3,1){\makebox(0,0){$d$}}
\put(1.3,0.5){\makebox(0,0){$G_{1}$}}
\multiput(7,0)(1,0){4}{\circle*{0.29}}
\multiput(8,1)(1,0){3}{\circle*{0.29}}
\put(7,0){\line(1,0){3}}
\put(8,1){\line(1,0){1}}
\put(8,0){\line(0,1){1}}
\put(8,1){\line(1,-1){1}}
\put(9,0){\line(1,1){1}}
\put(6.3,0.5){\makebox(0,0){$G_{2}$}}
\end{picture}
\caption{$G_{1}$ and $G_{2}$ are\ unicycle graphs, but only $\Psi(G_{2})$ is a
greedoid.}
\label{fig223}
\end{figure}
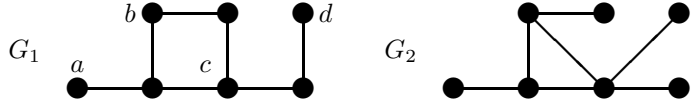

For instance, $S_{1}=\{a,b,d\}$ and $S_{2}=\{b,c,d\}$ belong to
$\Psi(G_{1}) $, where $G_{1}$ is the graph in Figure \ref{fig223}, but only
$S_{1} $\ has an accessibility chain, namely, $\{a\}\subset\{a,d\}\subset
S_{1}$. Nevertheless, having a simplicial vertex is a necessary but not a
sufficient condition for a stable set to admit an accessibility chain; e.g.,
$S_{2}=\{b,c,d\}\in\Psi(G_{1})$ has a pendant vertex and no accessibility
chain. However, there exist graphs where every maximum stable set has an
accessibility chain, e.g., the graph $G_{2}$ from Figure \ref{fig223}.

Evidently, if $\Psi(G)$ has the accessibility property, then every $S\in
\Psi(G)$, $\left\vert S\right\vert =k\geq2$, has an accessibility chain.

In this paper we characterize the unicycle graphs whose family of
local maximum stable sets are greedoids. Namely, we demonstrate that if
$C_{k}$ is the unique cycle of $G$, then the family $\Psi(G)$ is a greedoid
for $k=3$, while for $k\geq4$, $\Psi(G)$ is a greedoid if and only if either
\emph{(a)} $k$ is an even number and all maximum matchings of $G$ are uniquely
restricted, or \emph{(b)} $k$ is an odd number and the closed neighborhood of
every local maximum stable set of $G$ induces a K\"{o}nig-Egerv\'{a}ry graph.

\section{Results}

Let $C$ be the unique cycle of a graph $G$. Clearly, for every $e\in E(C)$,
the resulting graph $G-e$ is a forest.

\begin{theorem}
\label{th2}\cite{LevMan2} For any forest $T,\Psi(T)$ is a greedoid on its
vertex set.
\end{theorem}

\begin{figure}[h]
\setlength{\unitlength}{1.0cm} \begin{picture}(5,1.5)\thicklines
\multiput(5,0)(1,0){4}{\circle*{0.29}}
\multiput(6,1)(1,0){3}{\circle*{0.29}}
\put(5,0){\line(1,0){3}}
\put(6,1){\line(1,0){2}}
\put(7,0){\line(0,1){1}}
\put(8,0){\line(0,1){1}}
\put(4.7,0.3){\makebox(0,0){$a$}}
\put(5.7,0.3){\makebox(0,0){$b$}}
\put(7.2,0.3){\makebox(0,0){$d$}}
\put(5.7,1.3){\makebox(0,0){$c$}}
\put(7.2,1.3){\makebox(0,0){$e$}}
\put(8.2,0.3){\makebox(0,0){$f$}}
\put(8.2,1.3){\makebox(0,0){$g$}}
\end{picture}
\caption{$S=\{a,d,g\}\in\Psi(G)$, but $S$ admits no accessibility chain in $G
$.}%
\label{fig3535}%
\end{figure}
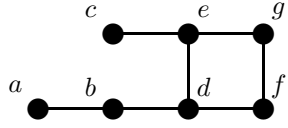

This assertion fails for general graphs, and even for unicycle graphs is not
always true (e.g., see the graph $G$ in Figure \ref{fig3535}, whose family
$\Psi(G)$ is not a greedoid).

\begin{proposition}
\label{lem55}If $\Psi(G)$ is a greedoid, then $\Omega(C_{k})\cap
\Psi(G)=\emptyset$ for its every induced cycle $C_{k}$ of size $k\geq4$.
\end{proposition}

\begin{proof}
Suppose that there exists $S\in\Omega(C_{k})\cap\Psi(G)$ in $G$ for some
$k\geq4$. Since $\Psi(G)$ is a greedoid, there is a chain of local maximum
stable sets
\[
\{x_{1}\}\subset\{x_{1},x_{2}\}\subset...\subset\{x_{1},x_{2},...,x_{k-1}%
\}\subset\{x_{1},x_{2},...,x_{q}\}=S,
\]
where $q=\left\vert S\right\vert \geq2$. Hence, $x_{1}$ must be a pendant
vertex in $G$, contradicting the fact that $x_{1}$ belongs to $V(C_{k})$.
\end{proof}

The graph $G$ from Figure \ref{fig3535} satisfies the condition that
$\Omega(C_{k})\cap\Psi(G)=\emptyset$ for its every cycle $C_{k}$ of size
$k\geq4$. Nevertheless, $\Psi(G)$ is not a greedoid, since $S=\{a,d,g\}\in\Psi(G)$, while
$S$ admits no accessibility chain. In other words, the converse of Proposition
\ref{lem55} is not true.

In the sequel, we distinguish between the following cases: $C=C_{3}$ and
$C=C_{k}$, $k\geq4$.

\begin{theorem}
\label{th5}If $G$ is a graph that has a $C_{3}$ as its unique cycle, then
$\Psi(G)$ is a greedoid.
\end{theorem}

\begin{proof}
Let $V(C_{3})=\{x_{i}:1\leq i\leq3\}$.

If $V(G)=V(C_{3})$, then it is easy to see that $\Psi(G)$ is a greedoid.

Let $V(G)\neq V(C_{3})$. According to Theorem \ref{th1}, it is sufficient to
show that $\Psi(G)$ satisfies the accessibility property. In other words, we
have to build an accessibility chain for any $S\in$ $\Psi(G)$.

Let $T_{ij},i\in\{1,2,3\},j\in\{0,1,...,n_{i}\}$ be subtrees of $G$ such that:
$T_{i0}=(\{x_{i}\},\emptyset)$, while for $j\geq1,T_{ij}$ is joined by an edge
to $x_{i},i\in\{1,2,3\}$, respectively (whenever such a subtree exists).

\begin{itemize}
\item \emph{Case 1.} $S\cap V(C_{3})=\emptyset$.
\end{itemize}

Then $S\in$ $\Psi(G-\{x_{1},x_{2},x_{3}\})$ and $T=G-\{x_{1},x_{2},x_{3}\}$ is
a forest. Therefore, by Theorem \ref{th2}, there is an accessibility chain for
$S$ in $T$. This is an accessibility chain for $S$ in $G$, as well, because
the neighborhoods of the sets belonging to the chain are the same in $T$ and
$G$.

\begin{itemize}
\item \emph{Case 2.} $S\cap V(C_{3})\neq\emptyset$, e.g. $S\cap V(C_{3}%
)=\{x_{1}\}$.
\end{itemize}

Let us denote
\begin{align*}
G_{1}  &  =G[\cup\{V(T_{1j}):0\leq j\leq n_{1}\}],\quad S_{1}=S\cap
V(G_{1}),\\
G_{2}  &  =G[\cup\{V(T_{2j}):1\leq j\leq n_{2}\}],\quad G_{3}=G[\cup
\{V(T_{3j}):1\leq j\leq n_{3}\}],\\
S_{ij}  &  =S\cap V(T_{ij}),1\leq j\leq n_{i},\quad S_{i}=\cup\{S_{ij},1\leq
j\leq n_{i}\},i=2,3.
\end{align*}

\begin{description}
\item \emph{Claim 1.} $S_{2},S_{3}$ are local maximum stable sets in
$G_{2},G_{3}$, respectively.
\end{description}

Otherwise, suppose that, for instance, $S_{2}\notin\Psi(G_{2})$. Then, there
is $W\in\Omega(N_{G_{2}}[S_{2}])$ with $\left\vert S_{2}\right\vert
<\left\vert W\right\vert $. Since $N_{G_{2}}[S_{2}]\subseteq N_{G}[S_{2}]$ and
$x_{1}\notin N_{G}[S_{2}]$, it follows that $S_{1}\cup W\cup S_{3}$ is a
stable set in $N_{G}[S]$, but larger than $S$, in contradiction to the choice
$S\in$ $\Psi(G)$.

\begin{description}
\item \emph{Claim 2.} $S_{1}\in\Psi(G_{1})$.
\end{description}

Otherwise, there must exist some stable set $W\subseteq N_{G_{1}}[S_{1}]$ with
$\left\vert S_{1}\right\vert <\left\vert W\right\vert $. Hence, $W$ is also
stable in $G$, and since $N_{G_{1}}[S_{1}]\subseteq N_{G}[S_{1}]$, it follows
that $W\cup S_{2}\cup S_{3}$ is a stable set in $N_{G}[S]$, but larger than
$S$, in contradiction to the choice $S\in$ $\Psi(G)$.

\begin{description}
\item \emph{Claim 3.} There is an accessibility chain of $S$ in $G$.
\end{description}

We distinguish between the following two cases.

\begin{description}
\item \emph{Case 3.1.} $S_{1}-\{x_{1}\}\in\Psi(G_{1})$.
\end{description}

Hence, we infer that $S_{1}-\{x_{1}\}\in\Psi(G_{1}-\{x_{1}\})$, which together
with \textit{Claim 1} imply
\[
(S_{1}-\{x_{1}\})\cup S_{2}\cup S_{3}\in\Psi(G-\{x_{1},x_{2},x_{3}\}).
\]
Therefore, according to Theorem \ref{th2}, there exists an accessibility chain
for the local maximum stable set $(S_{1}-\{x_{1}\})\cup S_{2}\cup S_{3}$ in
$T=G-\{x_{1},x_{2},x_{3}\}$, because $T$ is a forest. This is an accessibility
chain for $(S_{1}-\{x_{1}\})\cup S_{2}\cup S_{3}$ in $G$, as well, because the
neighborhoods of the sets belonging to the chain are the same in $T$ and $G$.
Clearly, this gives rise to an accessibility chain for $S_{1}\cup S_{2}\cup
S_{3}$.

\begin{description}
\item \emph{Case 3.2.} $S_{1}-\{x_{1}\}\notin\Psi(G_{1})$.
\end{description}

Since $S_{1}\in\Psi(G_{1})$ (by \textit{Claim 2}) and $\Psi(G_{1})$ is a
greedoid (by Theorem \ref{th2}), there is $v\in S_{1}$, such that
$S_{1}-\{v\}\in\Psi(G_{1})$.

We assert that $(S_{1}-\{v\})\cup S_{2}\cup S_{3}\in\Psi(G)$.

Otherwise, there is a stable set $A\subset N_{G}[(S_{1}-\{v\})\cup S_{2}\cup
S_{3}]$ with
\[
\left\vert A\right\vert >\left\vert S_{1}\right\vert -1+\left\vert
S_{2}\right\vert +\left\vert S_{3}\right\vert .
\]
Therefore, either $x_{2}\in A$ or $x_{3}\in A$. Without lack of generality
suppose that $x_{2}\in A$. Hence, there exists a stable set $B\subset
N_{G}[S_{2}]$ in $G_{2}$ such that $\left\vert B\right\vert =\left\vert
S_{2}\right\vert $ and $B\cup\{x_{2}\}$ is stable in $G$. Since $S_{1}%
-\{x_{1}\}\notin\Psi(G_{1})$, there exists $S_{0}\subset N[S_{1}-\{x_{1}\}]$,
such that $\left\vert S_{0}\right\vert >\left\vert S_{1}-\{x_{1}\}\right\vert
$. Clearly, $x_{1}\notin S_{0}\subset N[S_{1}]$. Consequently, $S_{0}\cup
B\cup\{x_{2}\}\cup S_{3}$ is a stable set in $N_{G}[S]$ of size greater than
$\left\vert S\right\vert $, and that contradicts the choice of $S\in$
$\Psi(G)$.

Thus, we obtain
\[
S,(S_{1}-\{v\})\cup S_{2}\cup S_{3}\in\Psi(G),
\]%
\[
x_{1}\in S_{1}-\{v\}\in\Psi(G_{1}),\quad S_{2}\in\Psi(G_{2}),\quad S_{3}%
\in\Psi(G_{3}),
\]
and%
\[
S=S_{1}\cup S_{2}\cup S_{3}\supset(S_{1}-\{v\})\cup S_{2}\cup S_{3}.
\]
Now, if $(S_{1}-\{v\})-\{x_{1}\}\in\Psi(G_{1})$, we continue as in
\textit{Case 1}, which leads immediately to an accessibility chain.

Otherwise, we find some vertex
\[
v^{\prime}\in S_{1}-\{v\}\in\Psi(G_{1}),v^{\prime}\neq x_{1},
\]
such that $(S_{1}-\{v\})-\{v^{\prime}\}$ belongs to $\Psi(G_{1})$, and we can
continue, as in \textit{Case 2}, to increase the size of the chain of local
maximum stable sets we are building one by one.

Since the set $S$ is finite, at the end of the above procedure we obtain an
accessibility chain of $S$ in $G$.
\end{proof}

Let us remark that the graph $G$ in Figure \ref{fig3535} is a unicycle
bipartite graph, and as we mentioned before, its family $\Psi(G)$ is not a
greedoid. However, there exist bipartite graphs whose families $\Psi(G)$ are
greedoids. 

Trying to characterize these bipartite graphs, we found out an
interesting connection between their local maximum stable sets of a graph and
their matchings, but of some special kind \cite{LevMan5,LevMan2,LevMan4,LevMan6}. 

A \textit{perfect matching} is a matching saturating all the vertices of the graph. 
A matching $M $ of a graph $G$ is called \textit{a uniquely restricted matching} if $M$
is the unique perfect matching of the subgraph induced by the vertices it
saturates \cite{GolHiLew}. 

Recall that a cycle $C$ is \textit{alternating} with respect to a matching $M$ if for any two incident edges of $C$
exactly one of them belongs to $M$, \cite{Krogdahl}. It is clear that an
$M$-alternating cycle should be of even size.

\begin{theorem}
\cite{GolHiLew}\label{th3} A matching $M$ in a graph $G$ is uniquely
restricted if and only if $G$ does not contain an alternating cycle with
respect to $M$.
\end{theorem}

Notice that the graph $G$ in Figure \ref{fig3535} has maximum matchings that
are not uniquely restricted. It turns out that the existence of such matchings is the real reason why
$\Psi(G)$ is not a greedoid.

\begin{theorem}
\label{th4}\cite{LevMan4} For a bipartite graph $G$, the family $\Psi(G)$ is a greedoid if
and only all maximum matchings of $G$ are uniquely restricted.
\end{theorem}

According to Theorem \ref{th3}, all maximum matchings of a unicycle
non-bipartite graph $G$ are uniquely restricted. 

\begin{figure}[h]
\setlength{\unitlength}{1.0cm} 
\begin{picture}(5,1.2)\thicklines
\multiput(2,0)(1,0){4}{\circle*{0.29}}
\multiput(3,1)(1,0){2}{\circle*{0.29}}
\put(2,0){\line(1,0){3}}
\put(3,0){\line(0,1){1}}
\put(3,1){\line(1,0){1}}
\put(4,1){\line(1,-1){1}}
\put(4.35,1){\makebox(0,0){$u$}}
\put(3.9,0.3){\makebox(0,0){$v$}}
\put(1,0.5){\makebox(0,0){$G_{1}$}}
\put(6.5,0.5){\makebox(0,0){$G_{2}$}}
\multiput(7.5,0)(1,0){4}{\circle*{0.29}}
\multiput(7.5,1)(1,0){4}{\circle*{0.29}}
\put(7.5,0){\line(1,0){3}}
\put(7.5,0){\line(1,1){1}}
\put(7.5,1){\line(1,0){3}}
\put(9.5,0){\line(0,1){1}}
\end{picture}
\caption{Non-bipartite triangle-free graphs {with unique perfect matchings}.}
\label{fig27}
\end{figure}
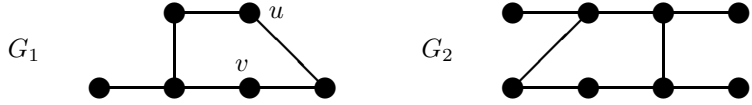

Nevertheless, it is not
sufficient for $\Psi(G)$ to be a greedoid. For example, the graphs
$G_{1},G_{2}$ in Figure \ref{fig27} are unicycle non-bipartite graphs,
$\Psi(G_{2})$ is a greedoid, while $\Psi(G_{1})$ is not a greedoid, because
$\{u,v\}\in\Psi(G_{1})$, but $\{u\},\{v\}\notin\Psi(G_{1})$.

A \textit{triangle-free graph} is a graph having no induced subgraph
isomorphic to $C_{3}$. The graphs in Figure \ref{fig27} are
triangle-free K\"{o}nig-Egerv\'{a}ry graphs, while $G_{1}[\{u,v\}]$ is not a
K\"{o}nig-Egerv\'{a}ry graph. The existence of such a subgraph is the reason for
$\Psi(G_{1})$ not to be a greedoid. 

\begin{theorem}
\cite{LevMan6}
\label{th7} If $G$ is a triangle-free graph, then $\Psi(G)$ is a greedoid if
and only if all its maximum matchings are uniquely restricted and the closed
neighborhood of each local maximum stable set of $G$ induces a
K\"{o}nig-Egerv\'{a}ry graph.
\end{theorem}

Now, using the fact that, by Theorem \ref{th3}, all the maximum matchings of a
graph without even cycles must be uniquely restricted, and combining Theorems
\ref{th5}, \ref{th4}, \ref{th7}, we conclude with the following.

\begin{theorem}
Let $C_{k}$ be the unique cycle of the graph $G$. Then, the following
statements are true:

\emph{(i)} if $k=3$, then $\Psi(G)$ is a greedoid;

\emph{(ii)} if $k=2q\geq4$, then $\Psi(G)$ is a greedoid if and only if all
maximum matchings of $G$ are uniquely restricted;

\emph{(iii)} if $k=2q+1\geq5$, then $\Psi(G)$ is a greedoid if and only if the
closed neighborhood of every local maximum stable set of $G$ induces a
K\"{o}nig-Egerv\'{a}ry graph.
\end{theorem}

\section{Conclusions}

In this paper we have completely characterized the unicycle graphs whose families of
local maximum stable sets form greedoids on their vertex sets. In
\cite{LevMan5} we showed that even unicycle graphs whose families of local
maximum stable sets are greedoids can be recognized in polynomial time. The
key question that remains open is whether there exists a polynomial time
recognition algorithm for odd unicycle graphs whose families of local maximum
stable sets are greedoids.

\end{document}